\documentclass[graybox]{svmult}

\usepackage{mathptmx}       
\usepackage{helvet}         
\usepackage{courier}        
\usepackage{type1cm}        
\usepackage{amssymb, amsmath, epsfig}
\usepackage{makeidx}         
\usepackage{graphicx}        
\usepackage{multicol}        
\usepackage[bottom]{footmisc}
\usepackage{algorithm,algorithmic}
\usepackage{todonotes}
\usepackage{hyperref}
\makeindex             
                       
 \newcommand{\R}{\mathbb R}
  \newcommand{\Gval}{{\sc gark}}
    \newcommand{\PGval}{{\sc pgark}}
  \newcommand{\BT}{{\sc bt}}
  \newcommand{\POD}{{\sc pod}}

\begin{document}
\bibliographystyle{siam}

\title*{Order reduction approaches for the algebraic Riccati equation and the LQR problem\thanks{Version of
Oct 20, 2017}}
\author{Alessandro Alla and Valeria Simoncini}
\institute{Alessandro Alla \at Department of Mathematics, PUC-Rio, Rio De Janeiro, Brazil.\\ \email{alla@mat.puc-rio.br}
\and Valeria Simoncini \at Alma Mater Studiorum - Universita' di Bologna, Italy and IMATI-CNR, Pavia, Italy.\\
\email{valeria.simoncini@unibo.it}} 
%
%
\maketitle

\abstract{We explore order reduction techniques for solving the algebraic Riccati equation (ARE),
and investigating the numerical solution of the linear-quadratic regulator problem (LQR). 
A classical approach is to build a surrogate low dimensional model of the 
dynamical system, for instance by means of balanced truncation, and then solve the corresponding ARE. 
Alternatively, iterative methods can be used to directly solve the ARE and use its approximate
solution to estimate quantities associated with the LQR.
We propose a class of Petrov-Galerkin strategies that simultaneously reduce the
dynamical system while approximately solving the ARE by projection.
This methodology significantly generalizes a recently developed Galerkin method by using a pair of
projection spaces, as it is often done in model order reduction of dynamical systems.
Numerical experiments illustrate the advantages of the
new class of methods over classical approaches when dealing with large matrices.}


\section{Introduction} \label{sec:1}

Optimal control problems for partial differential equations (PDEs) are an extremely important 
topic for many industrial applications in different fields,
 from  aerospace engineering to economics. 
The problem has been investigated with different strategies as 
open-loop (see e.g. \cite{HPUU09}) or closed-loop (see e.g. \cite{AFV17,FF14,GP11}). 

In this work we are interested in feedback control for linear dynamical systems and quadratic cost functional which is known as the Linear Quadratic Regulator (LQR) problem. 
Although most models are nonlinear, LQR is still a very interesting and powerful 
tool, for instance in the stabilization of nonlinear models under perturbations, where 
a control in feedback form can be employed.  

The computation of the optimal policy in LQR problems requires the solution of
an algebraic Riccati equation (ARE), a quadratic matrix equation with the dimension of the dynamical system. 
This is a major bottleneck in the numerical treatment of
 the optimal control problem, especially for high dimensional systems such as those
stemming from the discretization of a three-dimensional PDE.

Several powerful solution methods for the ARE have been developed throughout 
the years for small dynamical systems, based on spectral decompositions. 
The large scale case is far more challenging, as the whole spectral space of
the relevant matrices cannot be determined because of memory and computational
resource limitations. For these reasons, this algebraic problem is a very
active research topic, and major contributions have been given in the past
decade. 
Different approaches have been explored: variants of the Newton method have been
largely employed in the past \cite{K68},\cite{BS10}, while only more recently reduction type methods have emerged as
a feasible effective alternative; see, e.g., \cite{Benner.Saak.survey13},\cite{SSM14} and
references therein.  The recent work \cite{Simoncini2016} shows that a Galerkin class of
reduction methods for solving the ARE can be naturally set into the framework of
model order reduction for the original dynamical system.


As already mentioned, the LQR problem is more complicated when dealing with PDEs because 
its discretization leads to a very large system of ODEs and, as a consequence, 
the numerical solution of the ARE is computationally more demanding.
To significantly lower these computational costs and memory requirements,
model order reduction techniques can be employed.
Here, we distinguish between two different concepts of reduction approaches.
	
 A first methodology projects the dynamical system into a low dimensional system whose dimensions are
 much smaller than the original one; see, e.g., \cite{BCOW.17}.
Therefore, the corresponding reduced ARE is practical and feasible to compute on a standard computer.
The overall methodology  thus performs a first-reduce-then-solve strategy. 
This approach has been investigated with different model reduction techniques like Balanced Truncation (\BT) in e.g.\cite{A05}, Proper Orthogonal Decomposition (\POD, \cite{Sir87,Vol11}) in e.g. \cite{AK01,KS16} and via the interpolation of the rational functions, see \cite{Gugercin2008a},\cite{GB17}, \cite{bG15}. A different approach has been proposed in \cite{SH17} where the basis functions are computed from the solution of the high dimensional Riccati equation in a many query context. 
We note that basis generation in the context of 
model order reduction for optimal control problems is an active research topic 
(see e.g. \cite{AGH16,ASH17,KS16,KV08}). 
Furthermore, the computation of the basis functions is made by a 
Singular Value Decomposition (SVD) of the high dimensional data which can be very expensive. 
One way to overcome this issue was proposed in \cite{AK17} by means of randomized 
SVD which is a fast and accurate alternative to the SVD, and it is based on random samplings.

A second methodology follows a reduce-while-solve strategy. In this context,
recent developments aim at reducing the original problem by subspace projection,
and determining an approximate solution in a low dimensional approximation space. 
Proposed strategies either explicitly reduce the quadratic equation 
(see, e.g., \cite{SSM14} and references therein),
or approximately solve the associated invariant subspace problem 
(see, e.g., \cite{Benner.Bujanovic.16} and its references).
As already mentioned, these recently developed methods have shown to be effective alternatives to
classical variants of the Newton method, which require the solution of a
linear matrix equation at each nonlinear iteration; see, e.g.,
\cite{Benner.Li.Penzl.08} for a general description.

The aim of this paper is to discuss and compare the aforementioned
model order reduction methodologies 
for LQR problems. In particular, we compare the two approaches of 
reducing the dynamical system first versus building surrogate approximation 
of the ARE directly, using either Galerkin or Petrov-Galerkin projections.
The idea of using a Petrov-Galerkin method for the ARE appears to be
new, and naturally expands the use of two-bases type order reduction methods typically
employed for transfer function approximation.

To set the paper into perspective we start recalling LQR problem and its order reduction
in Section \ref{sec:2}. In Section \ref{sec:mordyn} we describe reduction strategies
of dynamical systems used in the small size case, such as proper orthogonal decomposition and
balanced truncation. Section \ref{sec:mor_ric} discusses the new class of projection strategies that
attack the Riccati equation, while delivering a reduced order model for the dynamical system.
Finally, numerical experiments are shown in Section \ref{sec:4} and conclusions are derived
 in Section \ref{sec:con}.

\section{The linear-quadratic regulator problem and model order reduction} \label{sec:2}
In this section we recall the mathematical formulation of the LQR problem. 
We refer the reader for instance to classical books such as
e.g. \cite{Lancaster.Rodman.95} for a comprehensive description of what follows. 
We consider a linear time invariant system of ordinary differential equations of dimension $n$:
\begin{align}
\dot{x}(t) &= A x(t) + B u(t), \quad x(0)=x_0, \quad t>0, \label{sys}\\
y(t) &= C x(t) + D u(t),\nonumber
\end{align}
{with $A\in\R^{n\times n}, B \in\R^{n\times m}, C \in \R^{p\times n}$ and $D \in\R^{p\times m}$. }
Usually, $x(t):[0,\infty]\rightarrow \R^n$ is called the state, $u(t):[0,\infty]\rightarrow \R^m$ the input or control and $y(t): [0,\infty]\rightarrow \R^p$ the output. 
{Furthermore, we assume that $A$ is passive. This may be viewed as a restrictive hypothesis, 
since the problems we consider only require that $(A, B)$ are stabilizable and $(C,A^T)$ controllable,
however this is convenient to ensure that the methods we analyze are well defined.}
In what follows, without loss of generality, we will consider $D\equiv 0$.
We also define the transfer function for later use:
\begin{equation}\label{transf}
G(s) = C(sI - A)^{-1} B.
\end{equation}
Next, we define the quadratic cost functional for an infinite horizon problem:
\begin{equation}\label{cost:qua}
J(u):= \int_0^\infty y(t)^T y(t) + u(t)^T R u(t)\, dt ,
\end{equation}
 where $R\in\R^{m\times m}$ is a symmetric positive definite matrix. The optimal control problem reads:
 \begin{equation}\label{prb_min}
\min_{u\in\R^m} J(u) \quad \text{such that} \quad x(t)\quad  \mbox{ solves }\quad  \eqref{sys}.
\end{equation}
The goal is to find a control policy in feedback form as:
\begin{equation}\label{opt_con}
u(t) = - K x(t) = - R^{-1} B^T P x(t)
\end{equation}
with the feedback gain matrix $K\in\R^{m\times n}$ and $P\in\R^{n\times n}$ is the unique symmetric and positive (semi-)definite 
matrix that solves the following ARE:
\begin{equation}\label{eq:ricc}
A^T P + P A - PBR^{-1}B^T P + C^T C =0,
\end{equation}
which is a quadratic matrix equation for the unknown $P$. 

We note that the numerical approximation of equation \eqref{eq:ricc} can be very expensive for large $n$. Therefore, we aim at the reduction of the numerical complexity by projection methods.

Let us consider a general class of tall matrices $V,W\in\R^{n\times r}$, whose 
columns span some approximation spaces. We chose these two matrices  such that they
are biorthogonal, that is $W^T V = I_r$.
Let us now assume that the matrix $P$, solution of \eqref{eq:ricc}, can be approximated as
$$
P\approx  W P_r W^T.
$$ 
Then the residual matrix can be defined as
$$
\mathcal{R}(P_r) = A^T W P_r W^T + W P_r W^T A - W P_r W^T  B R^{-1}B^T   W P_r W^T + C^TC .
$$
The small dimensional matrix $P_r$ can be determined by imposing the so-called
Petrov-Galerkin condition, that is orthogonality of the residual with respect
to range$(V)$, which in matrix terms can be stated as $V^T \mathcal{R}(P_r)V = 0$. Substituting
the residual matrix and exploiting the bi-orthogonality of $V$ and $W$ we obtain:

\begin{equation}\label{red:are}
A_r^T P_r + P_r A_r - P_r B_r R^{-1} B_r^T P_r + C_r^T  C_r =0
\end{equation}
%
where 
$$
A_r= W^T A V,\quad B_r = W^TB,\quad C_r = CV.
$$ 
It can be readily seen that equation \eqref{red:are} is
again a matrix Riccati equation, in the unknown matrix $P_r\in\R^{r\times r}$, of much smaller dimension
than $P$, provided that $V$ and $W$ generate small spaces. We refer to this equation as the
{\it reduced Riccati equation}. The computation of $P_r$ allows us to formally obtain
the approximate solution $W P_r W^T$ to the original Riccati equation \eqref{eq:ricc}, although
the actual product is never computed explicitly, as the approximation is kept in factorized
form.

The optimal control for the reduced problem reads 
$$u_r(t) = -K_r x_r(t) = - R^{-1}B_r^T P_r x_r(t)$$
with the reduced feedback gain matrix given by $K_r = KV \in \R^{m\times r}$.
Note that this $u_r(t)$ is different from the one obtained by first approximately
solving the Riccati equation and with the obtained matrix defining an approximation to $u(t)$;
see \cite{Simoncini2016} for a detailed discussion.

The Galerkin approach is obtained by choosing 
$V=W$ with orthonormal columns when imposing the condition on the residual.

To reduce the dimension of the dynamical system \eqref{sys}, we assume to approximate the full state vector as $x(t) \approx V x_r(t)$  with a basis matrix $V \in \R^{n\times r}$, where
$x_r(t):[0,\infty)\rightarrow \R^r$ are the reduced
coordinates. Plugging this ansatz into the dynamical
system~\eqref{sys}, and requiring a so called Petrov-Galerkin
condition yields

\begin{eqnarray} \label{eq:dynconautred}
  \dot{x}_r(t)&=& A_r x_r(t) + B_r u(t),\quad x_r(0)=W^Tx_0, \quad t>0,\\
  y_r(t) &=& C_r x_r.\nonumber
\end{eqnarray}

The reduced transfer function is then given by:
\begin{equation}\label{red:transf}
G_r(s) = C_r(sI_r - A_r)^{-1} B_r.
\end{equation}

The presented procedure is a generic framework for model reduction. It is
clear that the quality of the approximation depends on the
approximation properties of the reduced spaces. In the following sections, we will distinguish between the methods that directly compute $V,W$ upon the dynamical systems 
(see Section~\ref{sec:mordyn}) and those that readily reduce the ARE (see Section~\ref{sec:mor_ric}).
In particular, for each method we will discuss both Galerkin and Petrov-Galerkin projections,  to 
provide a complete overview of the methodology. The considered general Petrov-Galerkin approach for the
ARE appears to be new.



\section{Reduction of the dynamical system}\label{sec:mordyn}

In this section we recall two well-known techniques as \POD\ and \BT\ to compute the projectors $W,V$ starting from the dynamical systems. 


\subsection{Proper Orthogonal Decomposition}
\label{sec:mor_pod}
A common approach is based on the snapshot form of \POD\ proposed
in \cite{Sir87}, which works as
follows. We compute a set of snapshots $x(t_1),\dots,x(t_k)$
of the dynamical system \eqref{sys} corresponding to
a prescribed input $\bar u(t)$ and different time instances $t_1,\ldots,t_k$ 
and define the \POD\ ansatz of order $r$ for the state $x(t)$ by
\begin{equation}\label{pod_ans}
x(t)\approx\sum_{i=1}^r {(x_r)}_i(t)\psi_i,
\end{equation}
where the basis vectors $\{\psi_i\}_{i=1}^r$ are obtained from the
SVD of the snapshot matrix
$X=[x(t_1),\ldots,x(t_k)],$ i.e. $X=\Psi\Sigma \Gamma^T$, and the first $r$
columns of $\Psi=(\Psi_1,\dots,\Psi_n)$ form the \POD\ basis functions of rank $r$. Hence
we choose the basis vectors $V=W=(\Psi_1,\dots,\Psi_r)$ for the reduction
in~\eqref{eq:dynconautred}.  
This technique strongly relies on the choice of a given input $u$, whose optimal selection is
usually unknown. In this work, we decide to collect snapshots following the approach suggested in \cite{KS16} as considers a linearization of the ARE (which corresponds to a Lyapunov equation). Therefore, the snapshots are computed by the following equation:
\begin{equation}\label{adj}
  \dot{x}(t)=A^T x(t),\quad x(0)=c_i,\quad \text{for } i=1,\ldots p
 \end{equation}
where $c_i$ is the $i-$th column of the matrix $C$. The advantage of this
approach is that equation \eqref{adj} is able to capture the dynamics of the adjoint
equation which is directly related to the optimality
conditions, and we do not have to chose a reference input $\bar u(t).$ In order to obtain the \POD\ basis, one has to simulate the
high dimensional system and subsequently perform a SVD. As a consequence,
the computational cost may become prohibitive for large scale problems.
Algorithm \ref{pod_boris} summarizes the method.

\begin{algorithm}
\caption{\POD\ method to compute the reduced Riccati}
\label{pod_boris}
\begin{algorithmic}[1]
\REQUIRE $A, C, r$
\FOR{$i = 1,\ldots, p$}
\STATE Simulate system \eqref{adj} with initial condition $c_i$.
\STATE Build the snapshots matrix $X=[X, x_i(t_1),\ldots,x_i(t_k)]$
\ENDFOR
\STATE Compute the reduced SVD of $X=V\Sigma W^T$
\STATE Solve the reduced Riccati equation \eqref{red:are} for $P_r$.
\end{algorithmic}
\end{algorithm}

\subsection{Balanced truncation}
\label{sec:mor_bt}

The \BT\ method is a well-established reduced order modeling technique for linear time invariant systems \eqref{sys}. We refer to \cite{A05} for a complete
description of the topic. It is based on the solution of
the reachability Gramian $R$
and the observability Gramian $O$
which solve, respectively, the following Lyapunov equations
\begin{equation}\label{lyp_eq}
 AR + RA^T +  B B^T = 0,\quad  A^TO+ OA +  C^T C = 0.
 \end{equation}
We determine the Cholesky factorization of the Gramians
\begin{equation}\label{chol:gram}
R=\Phi\Phi^T\qquad O=\Upsilon\Upsilon^T.
\end{equation}
Then, we compute the reduced SVD of the Hankel operator $\Upsilon^T\Phi$ and set
$$W=\Upsilon U\Sigma^{1/2},\qquad V=\Upsilon V\Sigma^{1/2},$$
where $U, V\in\R^{n\times r}$ are the first $r$ columns of the left and right singular vectors of the Hankel operator and $\Sigma=\mbox{diag}(\sigma_1,\ldots,\sigma_r)$ matrix of the first $r$ singular values.

The idea of \BT\ is to neglect states that are
both, hard to reach and hard to observe. This is done by eliminating
states that correspond to low Hankel singular values $\sigma_i$.
This method is very popular in the small case regime, also because the whole procedure can 
be verified by a-priori error bounds in several system norms, and
the Lyapunov equations can be solved very efficiently. In the
large scale these equations need to be solved approximately;
see, e.g., ~\cite{Benner.Saak.survey13}.

In summary, the procedure first solves the two Lyapunov equations at
a given accuracy and then determines biorthogonal bases for the
reduction spaces by using a combined spectral decomposition of 
the obtained solution matrices. We note that this method is also very expensive for large $n$. The algorithm is summarized below in Algorithm \ref{alg:bt}.
\begin{algorithm}
\caption{\BT\ method to compute the reduced Riccati}
\label{alg:bt}
\begin{algorithmic}[1]
\REQUIRE $A,B,C$ and the dimension of the reduced problem $r$
\STATE Compute $R, O$ from \eqref{lyp_eq} and their Cholesky factorization \eqref{chol:gram}
\STATE Compute the reduced SVD of the Hankel operator
\STATE Set $W=\Upsilon U\Sigma^{1/2},\qquad V=\Upsilon V\Sigma^{1/2},$
\STATE Solve the reduced Riccati equation \eqref{red:are} for $P_r$.
\end{algorithmic}
\end{algorithm}

\section{Adaptive reduction of the algebraic Riccati equation}\label{sec:mor_ric}
In the previous section the reduced problem was obtained by a sequential procedure:
first system reduction of a fixed order $r$ and then solution of the reduced Riccati equation \eqref{red:are}.
A rather different strategy consists of determining the reduction
bases {\em while} solving the Riccati equation. In this way, we combine both the reduction
of the original system and of the ARE.
While the reduction bases $V$ and $W$ are being generated by means
of some iterative strategy, it is immediately possible
to obtain a reduced Riccati equation by projecting the problem onto
the current approximation spaces. The quality of the
two spaces can be monitored by checking how well the Riccati equation is solved by means of its residual; if 
the approximation is not satisfactory, the spaces can
be expanded and the approximation improved. 

The actual space dimensions are
not chosen a-priori, but tailored with the accuracy of the approximate
Riccati solution. Mimicking what is currently available in
the linear equation literature, the reduced problem can be obtained by imposing 
some constraints that uniquely identify an approximation. 
The idea is very natural and it was indeed presented in \cite{JaimoukhaFeb.1994},
where a standard Krylov basis was used as approximation space. However, only more
recently, with the use of rational Krylov bases,
a Galerkin approach has shown its real potential as a solver for the
Riccati equation; see, e.g., \cite{Heyouni.Jbilou.09},\cite{SSM14}. 
A more general Petrov-Galerkin approach was missing. We aim to fill this gap.
In the following we give more details on these procedures.

Given an approximate solution $\tilde P$ of \eqref{eq:ricc}
written as $\tilde P = W Y W^T$ for some $Y$ to be determined, the 
former consists to require that
the residual matrix is orthogonal to this same space, range($W$), so 
that in practice $V=W$. The Petrov-Galerkin procedure
imposes orthogonality with respect to the space range($V$), where $V$ is different
from $W$, but with the same number of columns.

In \cite{JaimoukhaFeb.1994} a first implementation of
 a Galerkin procedure was introduced, and the orthonormal columns
of $V$ spanning the (block) Krylov subspace 
${\cal K}_r(A^T, C^T)=
{\rm range}([ C^T, A^T C^T, \ldots, (A^T)^{r-1} C^T])$; see also
\cite{Jbilou_03} for a more detailed treatment and for numerical experiments.
Clearly, this definition generates a sequence of nested approximation spaces,
that is ${\cal K}_r(A^T, C^T)\subseteq {\cal K}_{r+1}(A^T, C^T)$,
whose dimension can be increased iteratively until a desired accuracy is achieved.
More recently, in \cite{Heyouni.Jbilou.09} and \cite{SSM14},
rational Krylov subspaces have been used, again in the Galerkin framework. In particular,
the special case of the {\it extended} Krylov subspace
${\cal K}_r(A^T, C^T) + {\cal K}_r((A^T)^{-1}, (A^T)^{-1}  C^T)$ was discussed
in \cite{Heyouni.Jbilou.09}, while the fully rational space
$$
{\cal K}_r(A^T, C^T, {\pmb\sigma}) := {\rm range}([C^T, (A^T-\sigma_2I)^{-1}C^T, \ldots,
 (A^T-\sigma_2I)^{-1} \cdots (A^T-\sigma_r I)^{-1}C^T])
$$
was used in \cite{SSM14}. The rational shift parameters ${\pmb\sigma}=
\{\sigma_2, \ldots, \sigma_r\}$ 
can be computed on the fly at low cost, by adapting the selection to the
current approximation quality \cite{Druskin.Simoncini.11}. 
Note that dim(${\cal K}_r(A^T, C^T, {\pmb\sigma}) ) \le r p$, where
$p$ is the number of columns of $C^T$.
In \cite{SSM14} it was also shown that a fully rational space can be more beneficial than the extended Krylov subspace for the Riccati equation.
In the following section
we are going to recall the general procedure associated
with the Galerkin approach, and introduce the algorithm for the Petrov-Galerkin
method, which to the best of the authors' knowledge is new.
In both cases we use the fully rational Krylov subspace with adaptive
choice of the shifts.  

It is important to realize that
${\cal K}_r(A^T, C^T, {\pmb\sigma})$ does not depend on the coefficient
matrix $BB^T$ of the second-order term in the Riccati equation.
Nonetheless, experimental evidence shows good performance. This issue was analyzed in \cite{Lin.Simoncini.15,Simoncini.16}
where however the use of the matrix $B$ during the computation of the parameters
was found to be particularly effective; a justification of this behavior
was given in \cite{Simoncini.16}. In the following we thus employ this 
last variant when using rational Krylov subspaces. More details will be given
in the next section.


\subsection{Galerkin and Petrov-Galerkin Riccati}\label{sec:petrgal_ricc}

In the Galerkin case, we will generate a matrix $W$ whose columns span the rational Krylov subspace
${\cal K}_r(A^T, C^T, {\pmb\sigma})$ in an iterative way, that is one block of columns at
the time. This can be obtained by an Arnoldi-type procedure; see, e.g., \cite{A05}.
The algorithm, hereafter {\sc gark} for Galerkin Adaptive Rational Krylov, works as follows:

\begin{algorithm}
\caption{\Gval\ method to compute the reduced Riccati equation}
\label{Gal-Ricc}
\begin{algorithmic}[1]
\REQUIRE $A,C, \pmb\sigma$
\FOR{$ r = 1,2,\ldots$}
\STATE Expand the space ${\cal K}_r(A^T,C^T, {\pmb\sigma})$;
\STATE Update the reduced matrices $A_r, B_r$ and $C_r$ with the newly
generated vectors;
\STATE Solve the reduced Riccati equation for $P_r$;
\STATE Check the norm of the residual matrix ${\cal R}(P_r)$
\STATE If satisfied stop with $P^{r}$ and the basis $W$ of ${\cal K}_r(A^T, C^T, {\pmb\sigma})$.
\ENDFOR
\end{algorithmic}
\end{algorithm}


The residual norm can be computed cheaply without the actual computation of the
residual matrix; see, e.g., \cite{SSM14}. The parameters $\sigma_j$
can be computed adaptively as the space grows; we refer the reader to
\cite{Druskin.Simoncini.11} and \cite{Simoncini.16} for more details.\\

In the general Petrov-Galerkin case, the matrix $W$ is generated the same way,
while we propose to compute
the columns of $V$ as the basis for the rational Krylov subspace
${\cal K}_r(A,B, {\pmb\sigma})$; note that the starting block is now $B$, and
the coefficient matrix is the transpose of the previous one.
The two spaces are now constructed and expanded at the same time, so that
the two bases can be enforced to be biorthogonal while they grow. For completeness,
we report the algorithm in the Petrov-Galerkin setting in Algorithm \ref{PG-Ricc}
(hereafter {\sc pgark} for Petrov-Galerkin Adaptive Rational Krylov).
\begin{algorithm}
\caption{\PGval\ method to compute the reduced Riccati equation}
\label{PG-Ricc}
\begin{algorithmic}[1]
\REQUIRE $A, B, C, \pmb\sigma$
\FOR{$ r = 1,2,\ldots$}
\STATE Expand the spaces ${\cal K}_r(A^T, C^T, {\pmb\sigma})$, ${\cal K}_r(A, B, {\pmb\sigma})$;
\STATE Update the reduced matrices $A_r, B_r$ and $C_r$ with the newly
generated vectors;
\STATE Solve the reduced Riccati equation for $P_r$;
\STATE Check the norm of the residual matrix ${\cal R}(P_r)$
\STATE If satisfied stop with $P^{r}$ and the basis $W$ of ${\cal K}_r(A^T, C^T, {\pmb\sigma})$.
\ENDFOR
\end{algorithmic}
\end{algorithm}

The parameters $\sigma_j$ are computed for one space and used also for the other space. In this more general case, the formula for the residual matrix norm is not as cheap
as for the Galerkin approach. 
We suggest the following procedure. We first recall that for
rational Krylov subspace ${\cal K}_r(A, B, {\pmb\sigma})$
the following relation holds (we assume here full dimension of the generated
space after $r$ iterations):
$$
A^T W = W A_r + \hat w a_r^T,  \quad
{a_r} \in \R^{(r+1) m },
$$
for certain vector $\hat w$ orthogonal to $W$, which changes as the iterations proceeds, that is
as the number of columns $W$ grows; we refer the reader to \cite{Lin.Simoncini.12} for a derivation
of this relation, which highlights that the distance of range($W$) from an invariant
subspace of $A$ is measured in terms of a rank-one matrix.
We write
{\footnotesize
\begin{eqnarray*}
\mathcal{R}(P_r) &=& A^T W P_r W^T + W P_r W^T A - W P_r W^T  B R^{-1}B^T   W P_r W^T + C^T C \\
&=& W A_r P_r  W^T + \hat w a_r^T P_r W^T + W P_r A_r^T W^T + W P_r a_r \hat w^T  
- W P_r B_r R^{-1} B_r^T P_r W^T + W E E^T W^T \\
&=& \hat w a_r^T P_r W^T + W P_r a_r \hat w^T \\
&=&
[W, \hat w] \begin{bmatrix} 0 & P_r a_r \\ a_r^T P_r  & 0 \end{bmatrix} [W, \hat w]^T, 
\end{eqnarray*}
}%
where we also used the fact that $C^T = W E$ for some matrix $E$.
If $[W, \hat w]$ had orthonormal columns, as is the case for Galerkin, then
$\|{\mathcal R}(P_r)\|^2 = 2 \|P_r a_r\|^2$, which can be cheaply computed.

To overcome the nonorthogonality of $[W, \hat w]$,
we suggest to perform a reduced QR factorization of $[W, \hat w]$ that maintains its columns orthogonal. This QR factorization does not have to 
be redone from scratch at each iteration, but it can be updated as the matrix $W$ grows.
If $[W, \hat w] = Q_W R_W$ with $R_W\in\R^{(r+1)m\times (r+1)m}$ upper triangular, then
$$
\|{\cal R}(P_r)\| =
\|R_W \begin{bmatrix} 0 & P_r a_r \\ a_r^T P_r  & 0 \end{bmatrix} R_W^T\|.
$$
The use of coupled (bi-orthogonal)
bases has the recognized advantage of explicitly using
both matrices $C$ and $B$ in the construction of the reduced
spaces. This coupled basis approach has been largely exploited in the approximation
of the dynamical system transfer function by solving
a multipoint interpolation problem; see, e.g., \cite{A05} for
a general treatment and \cite{Barkouki2015} for a recent implementation.
In addition, coupled bases can be used to simultaneously approximate
both system Gramians leading to a large-scale \BT\
strategy; see, e.g.,  \cite{Jaimoukha.Kasenally.95} for 
early contributions using bi-orthogonal standard Krylov subspaces\footnote{We are unaware
of any available implementation of rational Krylov subspace based approaches for 
large scale \BT\ either with single or coupled bases, that simultaneously performs
the balanced truncation while approximating the Gramians.}.
On the other hand, a Petrov-Galerkin procedure has several drawbacks
associated with the construction of the two bases. More precisely, the
two bi-orthogonal bases are generated by means of a 
Lanczos-type recurrence, which is known to have both stability
and breakdown problems in other contexts such as linear system and
eigenvalue solving. At any iteration it may happen that the
new basis vectors $w_j$ and $v_j$ are actually orthogonal or
quasi-orthogonal to each other, giving rise to a possibly incurable
breakdown \cite{GutknechtApr.1992}. We have occasionally
experienced this problem in our numerical tests, and it certainly occurs
whenever $CB=0$. 
In our specific context, an additional difficulty arises.
The projected matrix $A_r =W^T A^T V$ is associated with a
bilinear rather than a linear form, so that its
field of values may be unrelated to that of $A^T$. As a consequence,
it is not clear the type of hypotheses we need to impose
on the data to ensure that the reduced Riccati equation \eqref{red:are}
admits a unique stabilizable solution. 
Even in the case of $A$ symmetric, the two bases will be
different as long as $C \ne B^T$.
All these questions are
crucial for the robustness of the procedure
and deserve a more throughout analysis which will be the
topic of future research.

From an energy-saving standpoint, it is worth remarking that the
Petrov-Galerkin approach uses twice as many memory allocations than
the Galerkin approach, while performing about twice the number of
floating point operations.
In particular, constructing the two bases requires two
system solves, with $A^T-s_j I$ and with $A-\bar s_jI$
respectively, at each iteration. Therefore, unless convergence
is considerably faster, the Petrov-Galerkin approach may not be
superior to the Galerkin method in the solution of the ARE.

\section{Numerical Experiments}\label{sec:4}

In this section we present and discuss our numerical tests. 
We consider the discretization of the following linear PDE
\begin{equation}
\label{state:eq}
\begin{aligned}
w_t-\varepsilon \Delta w+ 
 \gamma w_x + \gamma w_y  &=\mathbf{1}_{\Omega_B} u&&\text{in } \Omega\times(0,+\infty),\\
w(\cdot,0)&=w_0&&\text{in } \Omega,\\
w(\cdot,t)&=0&&\text{in } \partial\Omega\times(0,\infty),
\end{aligned}
\end{equation}
where $\Omega\subset\mathbb R^2$ is an open interval, $w:\Omega\times[0, \infty]\rightarrow \R^2$ denotes the state, and the parameters $\varepsilon$ and $\gamma$ are real positive constants. 
The initial value is $w_0$ and the function $\mathbf{1}_{\Omega_B}$ is the indicator function over the domain $\Omega_B\subset\R^2$. Note that we deal with zero Dirichlet boundary conditions. 
The problem in \eqref{state:eq} includes the heat equation, 
for $\varepsilon \neq 0, \gamma=0$, and a class of convection-diffusion equations 
for $\varepsilon\neq 0$ and $\gamma \neq 0$.  
{Furthermore, we define an output of interest by: 
\begin{equation}\label{out:int}
 s(t):=\dfrac{1}{|\Omega_C|}\int_{\Omega_C} w(x,t)\, dt,
 \end{equation}
where $\Omega_C\subset\R^2.$ 
Space discretization of equation \eqref{state:eq} by standard centered finite differences together with a rectangular quadrature rule for \eqref{out:int} lead to a system of the form \eqref{sys}.}
%
%
In general, the dimension $n$ of the dynamical system \eqref{sys} is rather large (i.e., $n\gg1000$) 
and the numerical treatment of the
corresponding ARE is computationally expensive or even unfeasible. Therefore, 
model order reduction is appropriate to lower the dimension of the optimal control problem \eqref{prb_min}.
We will report experiments with small size problems, where all discussed methods can be
employed, and with large size problems, where only the Krylov subspace based strategies are applied.

The numerical simulations reported in this paper were performed on a Mac-Book Pro with 
1 CPU Intel Core i5 2.3 Ghz and 8GB RAM and the codes are written in MatlabR2013a.
In all our experiments, small dimensional Lyapunov and Riccati equations are solved
by means of built-in functions of the Matlab Control Toolbox.

Whenever appropriate, the quality of the current approximation of the ARE is monitored by using
the relative residual norm: 
\begin{equation}\label{res_ric}
{\tt R}_{P} = \dfrac{\|\mathcal{R}(P_r)\|_F}{\|C\|_F^2},
\end{equation}
and the dimension of the surrogate model $r$ is chosen such that ${\tt R}_P<10^{-8}$.
After the ARE solution is approximated the ultimate goal is to compute the feedback control \eqref{opt_con}. 
Therefore, we also report the error in the computation of feedback gain matrix $K$ as the iterations
proceed:

 \begin{equation}\label{err_K}
 \mathcal{E}_{K} = \dfrac{\|K_r - K \|_F}{\|K\|_F},
 \end{equation}
and we measure the quality of our surrogate model also by the $\mathcal{H}_2-$ error
 \begin{equation}\label{err_G}
 \mathcal{E}_{G} = \dfrac{\|G_r - G \|_{\mathcal{H}_2}}{\|G\|_{\mathcal{H}_2}}.
 \end{equation}
where
 $$
\|G(s)\|_{\mathcal{H}_2} := \dfrac{1}{2\pi} \left(\int_{-\infty}^{+\infty} \|G(i\omega)\|_F^2 d\omega\right)^{1/2}. 
$$
In particular, the approximation of the transfer function is one of the main targets of
reduced order modeling, where the reduced system is used for analysis purposes, while
the approximation of the feedback gain matrix is monitored to obtain a good control function.

\subsection{Test 1: 2D Linear heat equation}

In the first example we consider the linear heat equation. In \eqref{state:eq} we chose 
$\gamma= 0, \varepsilon =1, \Omega=[0,1]\times [0,1],$ and $\Omega_B=[0.2, 0.8] \times [0.2, 0.8]$.
In \eqref{sys}, the matrix $A$ is obtained by centered five points discretization. 
We consider a small problem stemming from
a spatial discretization step $\Delta x =0.05$, leading to a system of 
dimension $n= 441$. The matrix $C$ in \eqref{sys} is given by the indicator function over the domain 
$\Omega_C = [0.1, 0.9] \times [0.1, 0.9]$ and  $R\equiv I_m$ in \eqref{cost:qua}.

\begin{figure}[htbp]
\includegraphics[scale=0.2]{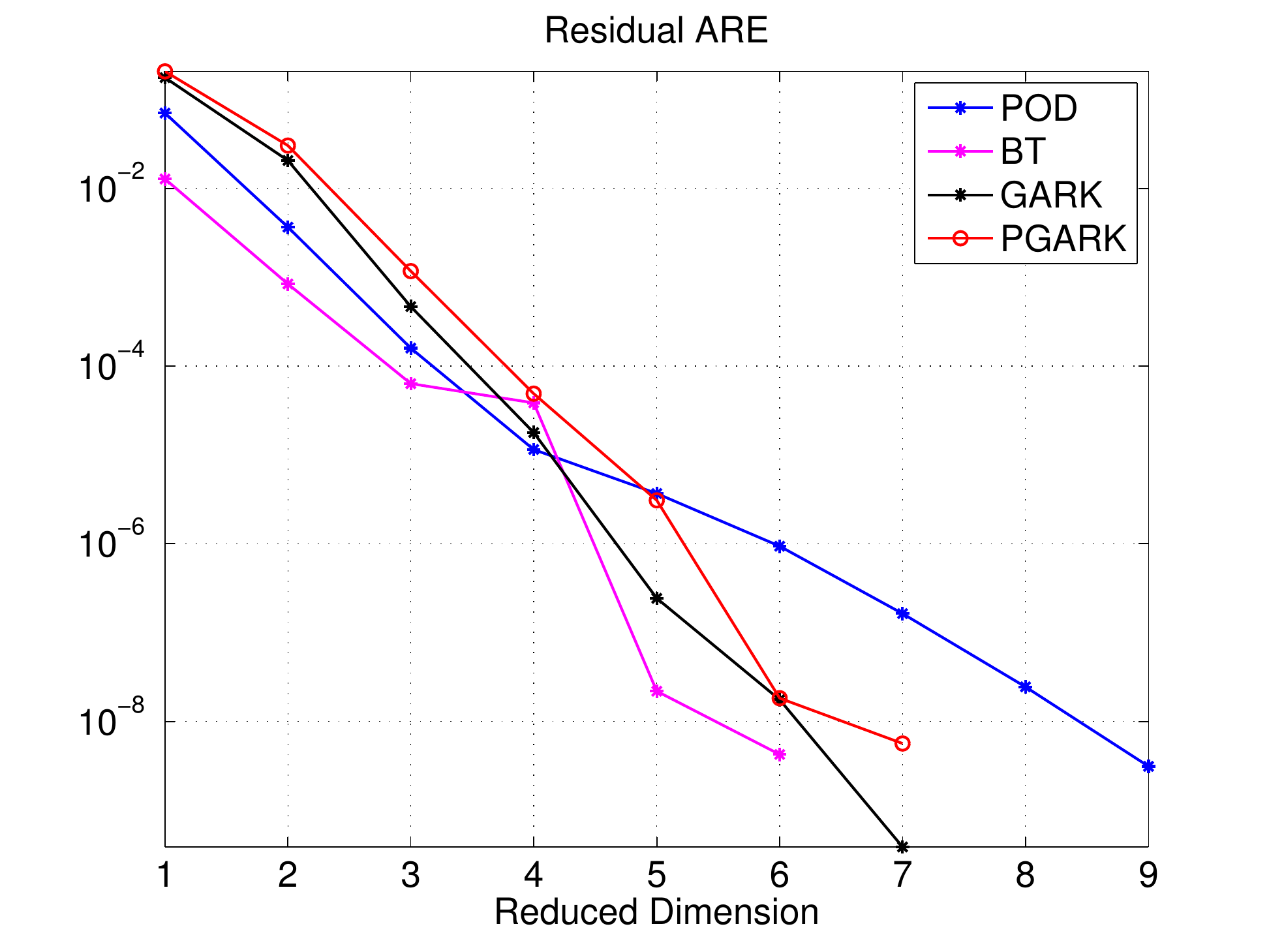}\includegraphics[scale=0.2]{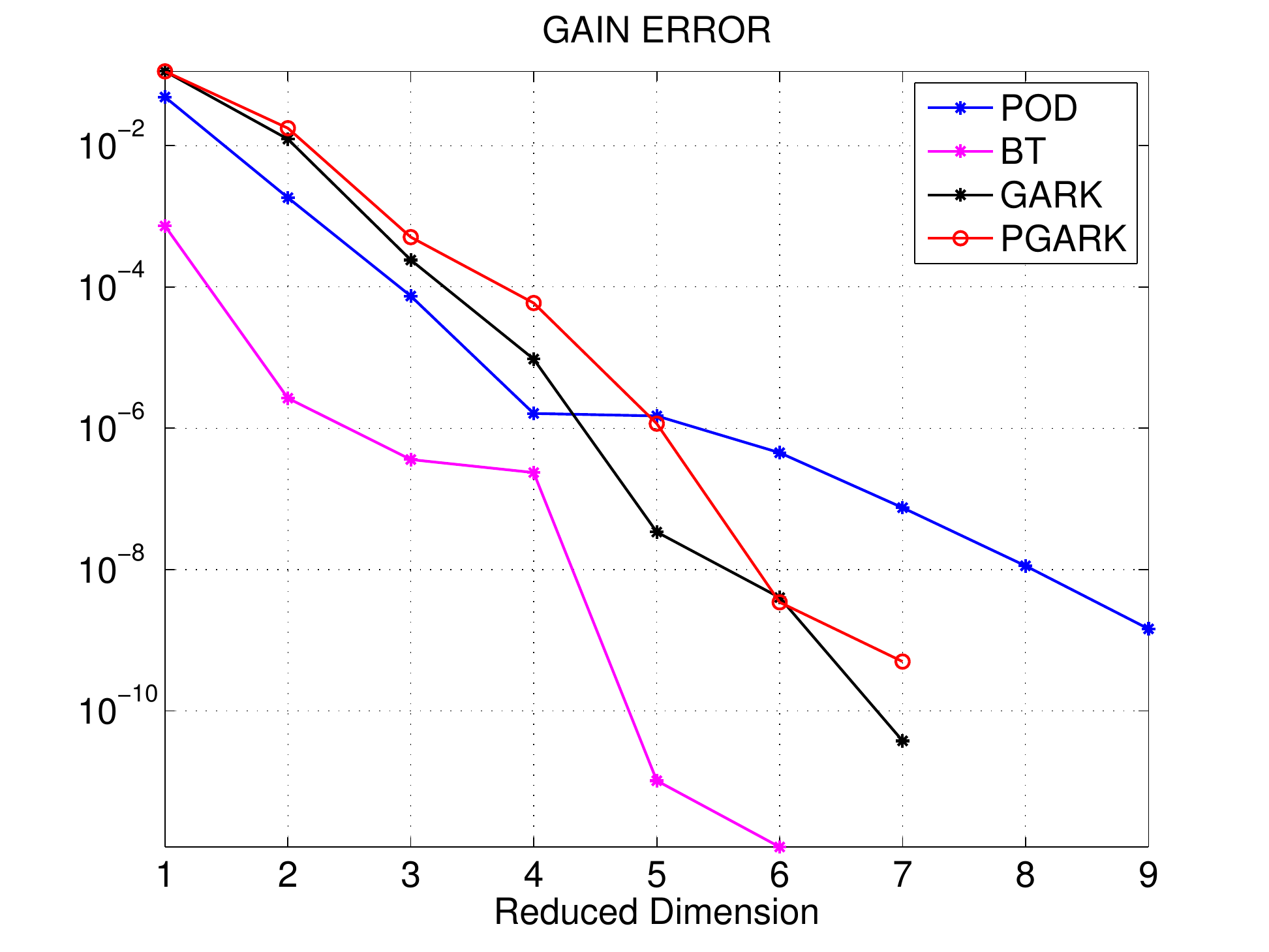}\includegraphics[scale=0.2]{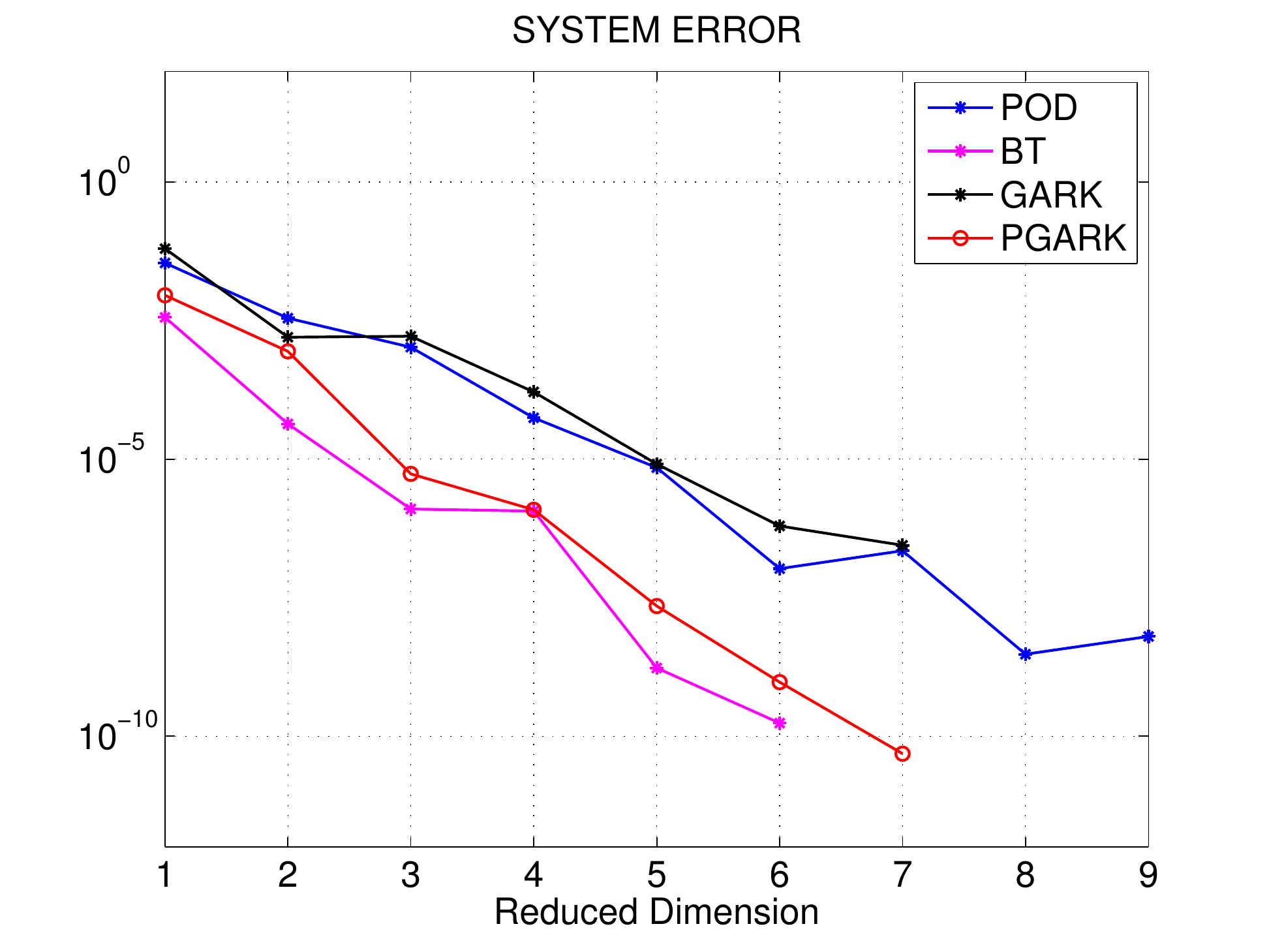}
\caption{Test 1: Convergence history of the relative residual norm ${\tt R}_P$ (left), Error ${\mathcal E}_K$ of the feedback gain matrix (middle), Error $\mathcal{E}_G$ for the approximation of reduced transfer function (right)}
\label{fig1}
\end{figure}

The left panel of Figure \ref{fig1} shows the residual norm 
history of the reduced ARE \eqref{red:are} when the
two projection matrices $V,W$ are computed by each of the four algorithms explained in the previous sections.  We can thus appreciate how the approximation proceeds as the reduced space is enlarged.
We note that \POD\ requires more basis functions to achieve the desired tolerance, 
whereas the \BT\ algorithm reaches it faster. However, the \POD\ method is a snapshots dependent method that is really influenced by the choice of the initial input $u(t)$ and the results may be different for other choices of the snapshots set. All the other proposed techniques are, on the contrary, input/output independent.

%
In the middle panel of Figure \ref{fig1}, we show how well the feedback gain matrix $K$ can be approximated with order reduction methods. It is interesting to see how the basis functions computed by \Gval\ and \PGval\ are able to approximate the matrix $K$ very well.  
%

Finally, we would like to show the quality of the computed basis 
functions in the approximation of the dynamical system in the right panel of Figure \ref{fig1}. 
In this example, \BT\ approximates the transfer function very well but this method requires the full accurate solution of two Lyapunov equations to be able to generate the reduced transfer function. In the large scale case this is clearly unfeasible. 

\begin{figure}[htbp]
\centering
\includegraphics[scale=0.3]{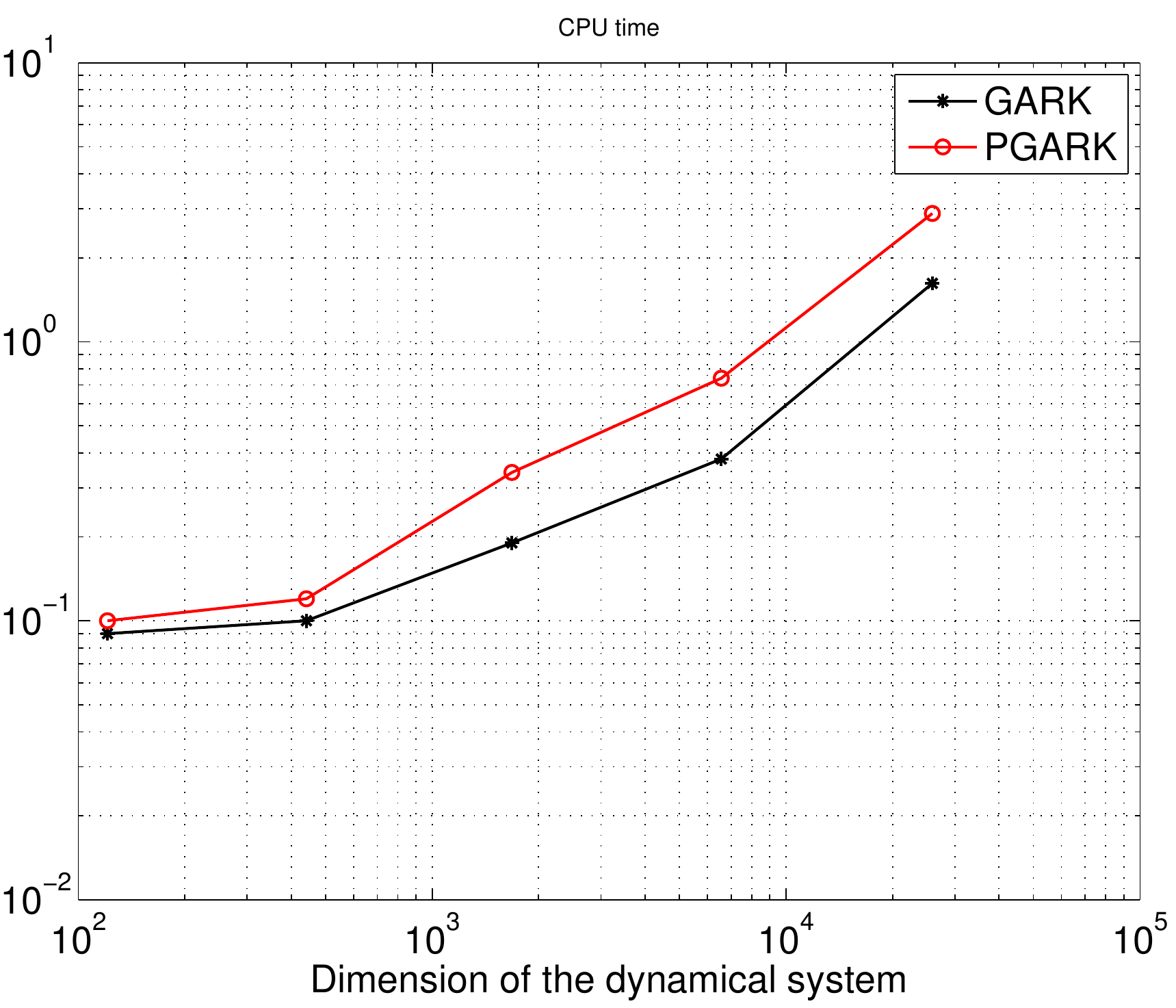}\includegraphics[scale=0.3]{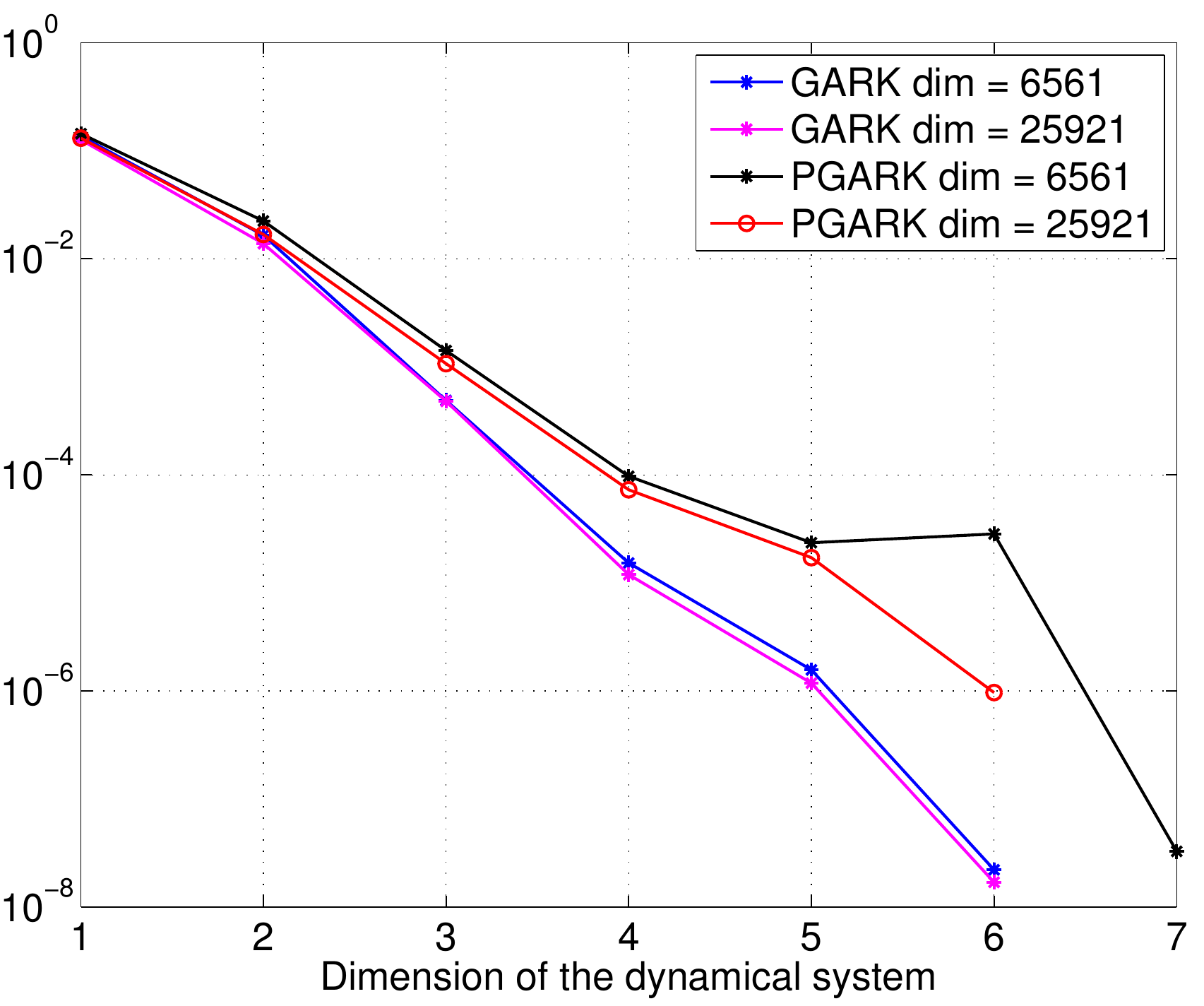}
\caption{Test 1: Comparison of \Gval\ and \PGval. Left: CPU time as the problem dimension $n$ increases. Right: relative residual norm history. }
\label{fig1:cpu}
\end{figure}

Last remark goes to the iterative methods \Gval\ and \PGval . We showed that, although the basis functions are built upon information of the ARE, they are also able to approximate the dynamical systems and the feedback gain matrix. This is a crucial point that motivates us to further investigate these methods in the context of the LQR problem. Furthermore, they are absolute feasible when the dimension $n$ of the dynamical system increases as shown in the left panel of Figure \ref{fig1:cpu}. We computed the CPU time of the iterative methods varying $\Delta x =\{0.1, 0.05, 0.025, 0.0125, 0,00625\}$. We note that both methods reach the desired accuracy in a few seconds even for $n = O(10^4)$. 
On the contrary both \POD\ and \BT\ would be way more expensive since their cost heavily depends on
the original dimension of the problem $n$. The right panel of 
Figure \ref{fig1:cpu} reports the history of the relative residual norm as the iterations proceed.


%




\subsection{Test 2: 2D Linear Convection-Diffusion Equation}
We consider the linear convection-diffusion equation in \eqref{state:eq} with 
$\gamma= 50, \varepsilon = 1, \Omega=[0,2]\times [0,2] \mbox{ and }\Omega_B=[0.2, 0.8] \times [0.2, 0.8].$
In \eqref{sys}, the matrix $A$ is given by centered five points finite difference
discretization plus an upwind approximation of the convection term (see e.g. \cite{Q12}). 
The spatial discretization step is $\Delta x =0.1$ and leads to a system of dimension $n= 441$. 
The matrix $C$ in \eqref{sys} is given by the indicator function 
over the domain $\Omega_C = [0.1, 0.9] \times [0.1, 0.9],$ and $R\equiv I_m$ in \eqref{cost:qua}.

\begin{figure}[htbp]
\includegraphics[scale=0.2]{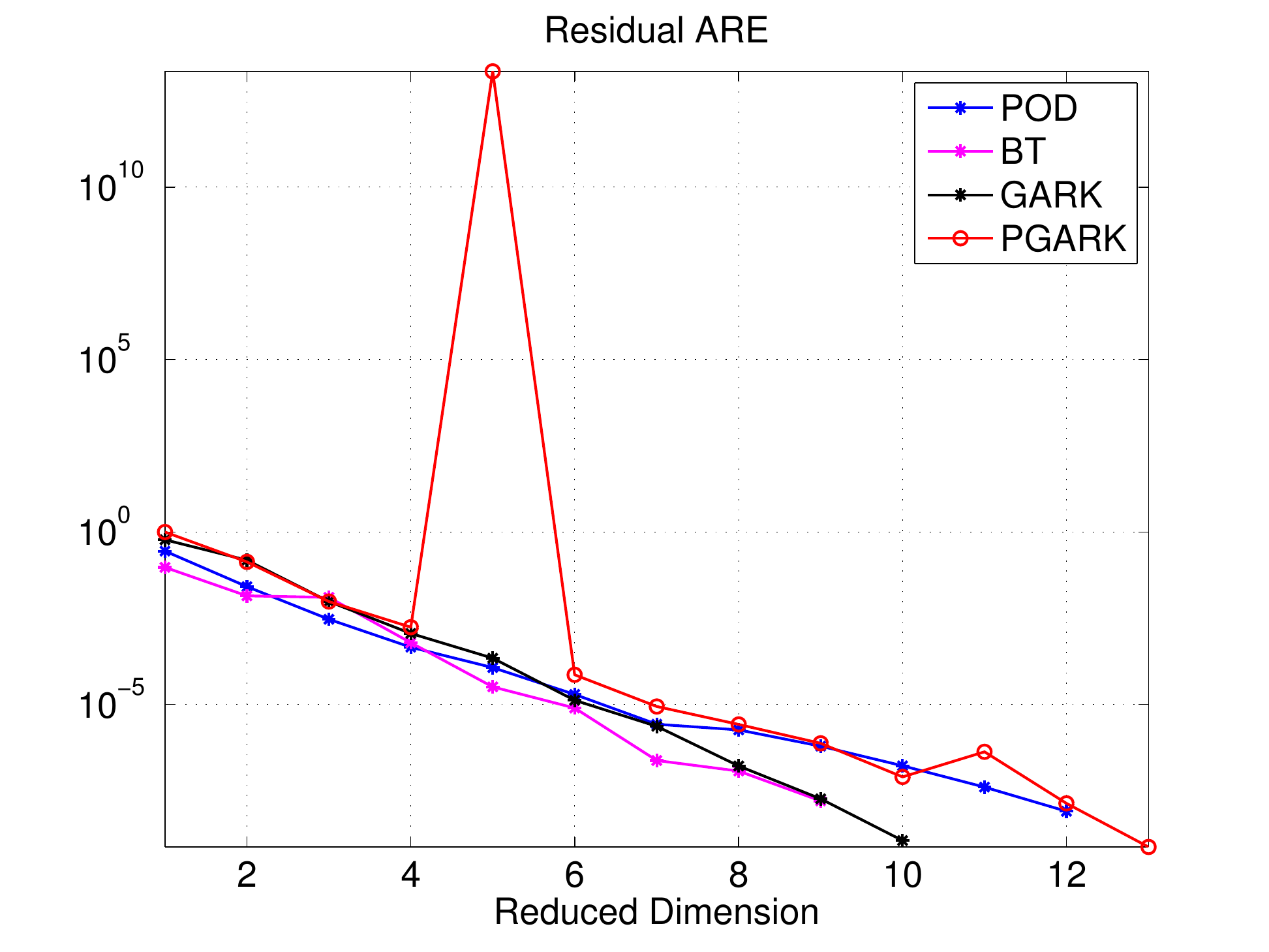}\includegraphics[scale=0.2]{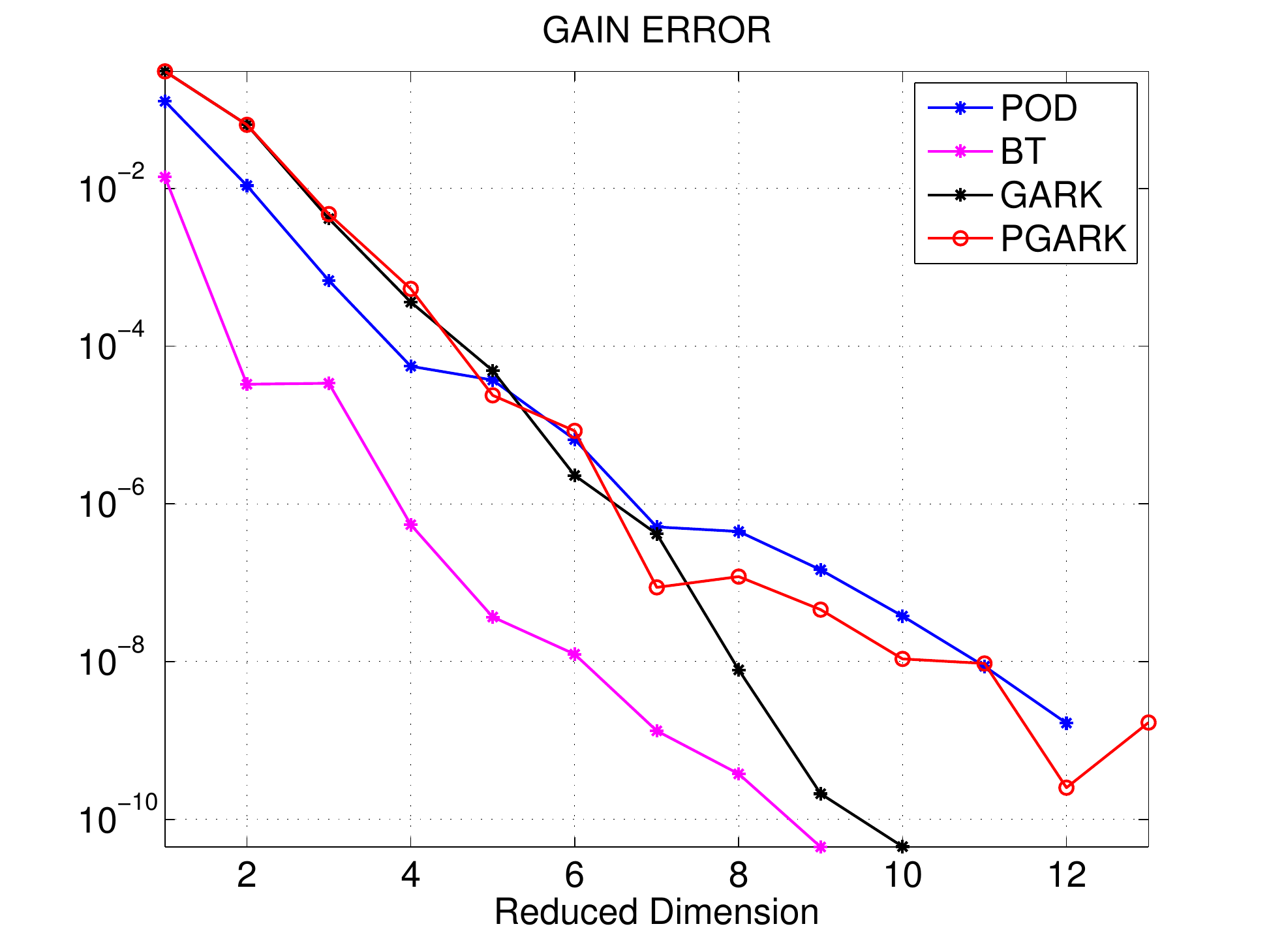}\includegraphics[scale=0.2]{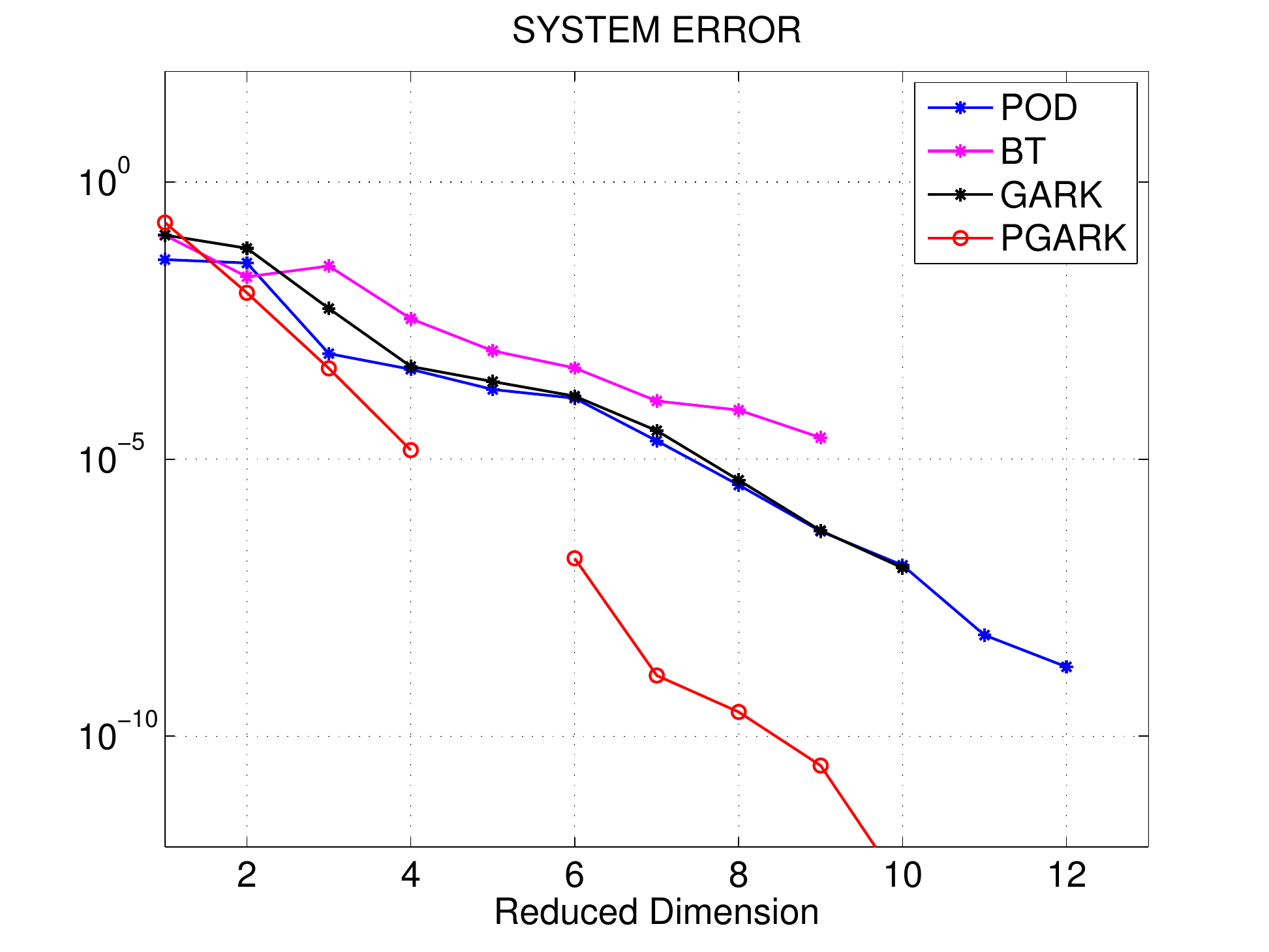}
\caption{Test 2: History of the relative residual norm (left), Error $\mathcal{E}_K$ of the feedback gain matrix (middle), Transfer function error $\mathcal{E}_G$ as the approximation space grows (right)}
\label{fig2}
\end{figure}

The left panel of Figure~\ref{fig2} shows the residual norm
history associated with the reduced ARE \eqref{red:are}, with different projection techniques.  
We note that all the methods converge with a different number of basis functions. We also note that in this example we can observe instability of the \PGval\ method as discussed in Section \ref{sec:petrgal_ricc}.

Middle panel of Figure \ref{fig2} reports the error in the approximation of the feedback gain matrix $K$.
It is very interesting to observe that even on this convection dominated problem all the methods can reach an accuracy of order $10^{-8}$.


Finally, we show the error in the approximation of the transfer function in the right panel of Figure \ref{fig2}. In this example, we can see that the 
\PGval\ method performs better than the others but it is rather unstable; this well known instability
problem will be analyzed in future work.
The discussion upon the quality of the basis function we had in Test 1 still hold true. The iterative methods are definitely a feasible alternative to well-known techniques as \BT\ and \POD.
%

\begin{figure}[htbp]
\centering
\includegraphics[scale=0.3]{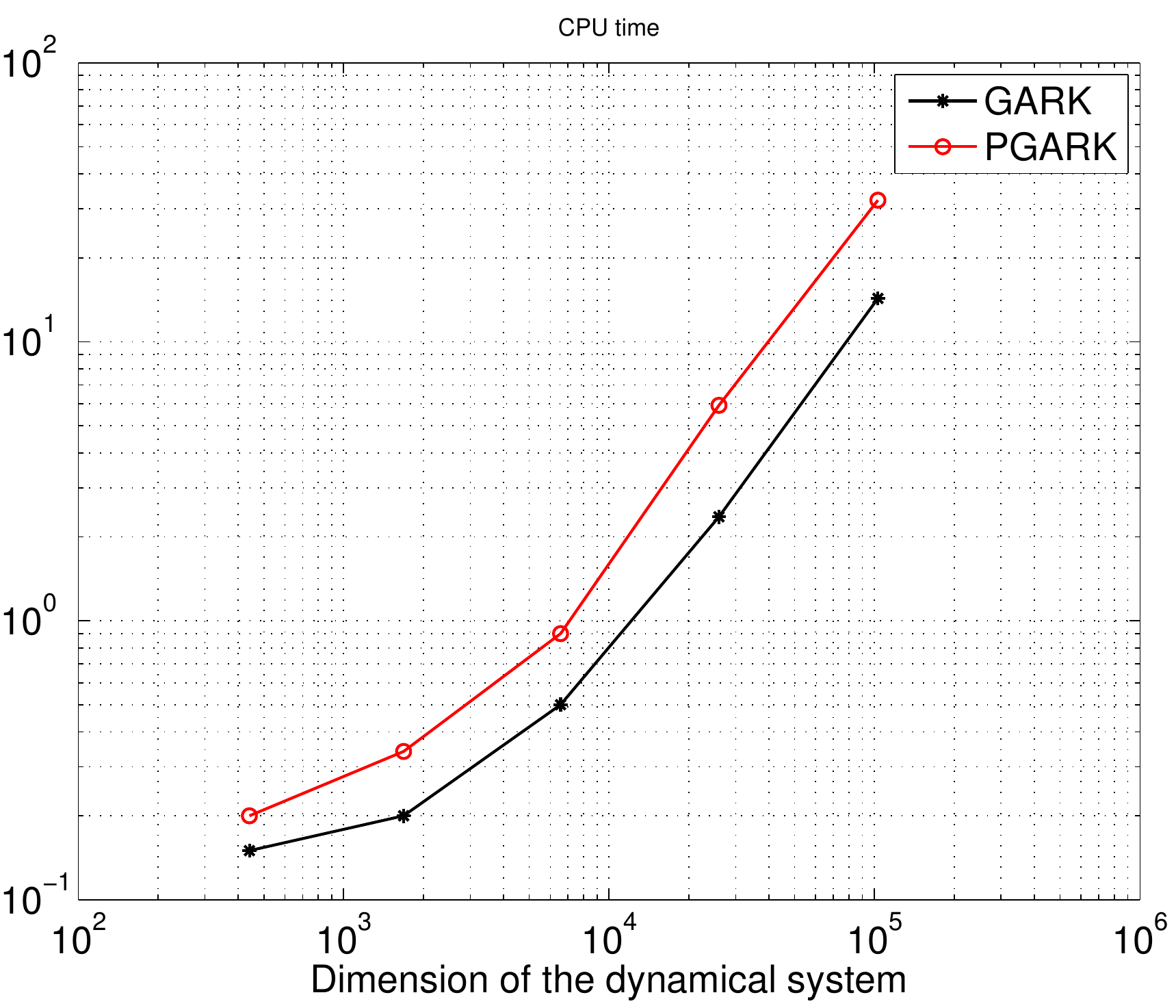}\includegraphics[scale=0.3]{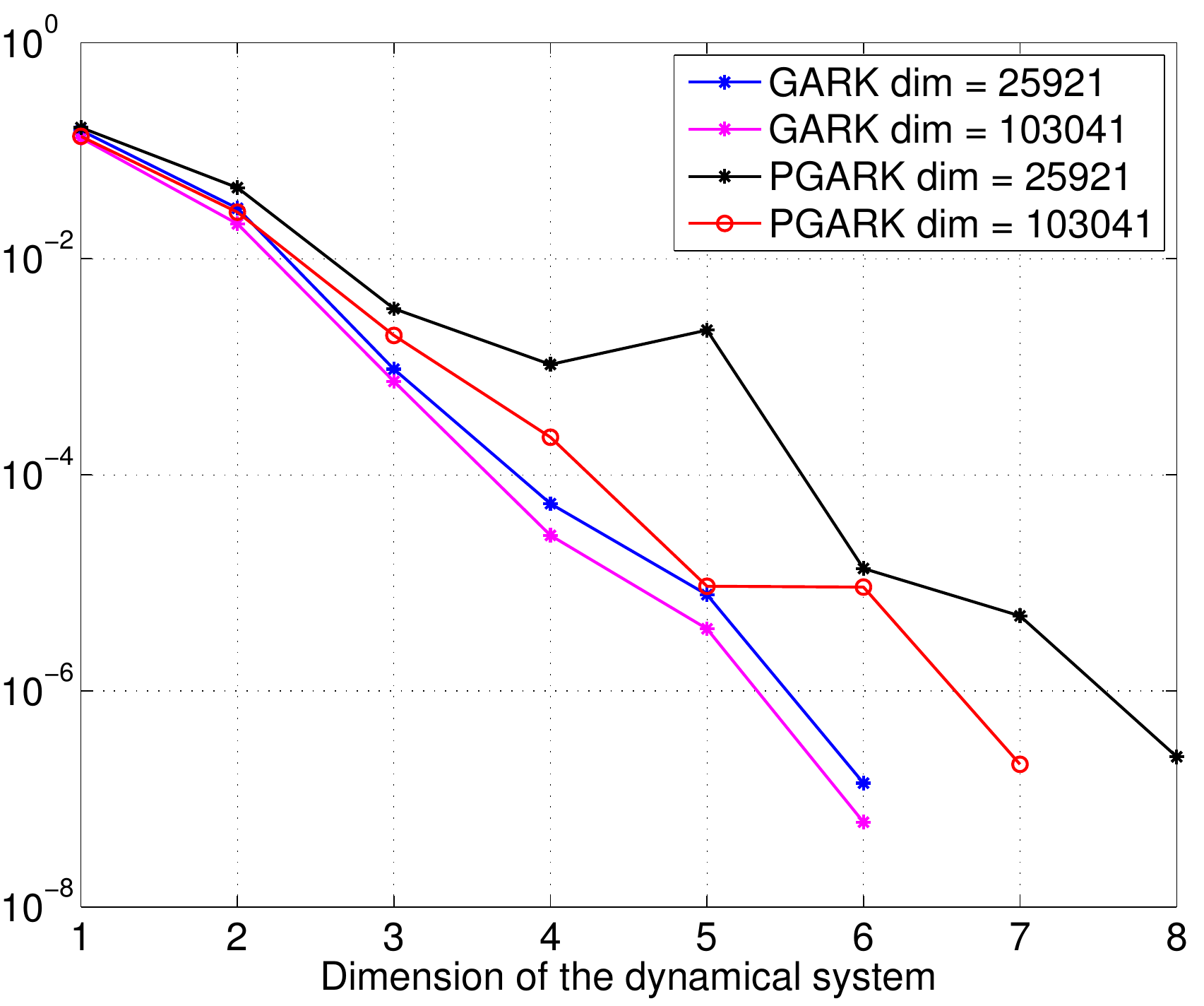}
\caption{Test 2: Comparison of \Gval\ and \PGval. Left: CPU time as the problem dimension $n$ increases.
Right: relative residual norm history.}
\label{fig2:cpu}
\end{figure}

In the left panel of Figure \ref{fig2:cpu}, we show the CPU time of \Gval\ and \PGval\
for different dimensions $n$ of the dynamical systems.
 We note that, again, the Galerkin projection reaches the desired accuracy faster.  
The right panel of Figure \ref{fig2:cpu} reports the history of the residual for both iterative methods, for
two different values of $n$.

%

\section{Conclusions}\label{sec:con}

We have proposed a comparison of different model order reduction techniques
for the ARE. We distinguished between two different strategies: 
(i) First reduction of the dynamical system complexity and then the
solution of the corresponding reduced ARE whereas;
 (ii) Simultaneous solution of the ARE and determination of the reduction spaces. 
The strength of the second strategy is its reliability even for very high dimensional problem,
where problems in the class (i) may have memory and computational difficulties. 
Experiments on both small and large dimensional problems confirm the promisingly
good approximation properties of rational Krylov methods as compared to more standard
approaches for the approximation of more challenging quantities in the LQR context.

\bibliography{%
Biblioteca, 
biblio_MF}

\begin{thebibliography}{10}

\bibitem{AFV17}
{\sc A.~Alla, M.~Falcone, and S.~Volkwein}, {\em Error analysis for {POD}
  approximations of infinite horizon problems via the dynamic programming
  approach}, SIAM Journal on Control and Optimization, 55 (2017),
  pp.~3091--3115.

\bibitem{AGH16}
{\sc A.~Alla, C.~Graessle, and M.~Hinze}, {\em A posteriori snapshot location
  for pod in optimal control of linear parabolic equations}, submitted,
  https://arxiv.org/abs/1608.08665,  (2016).

\bibitem{AK17}
{\sc A.~Alla and J.~Kutz}, {\em Randomized model order reduction}, submitted,
  https://arxiv.org/pdf/1611.02316.pdf,  (2017).

\bibitem{ASH17}
{\sc A.~Alla, A.~Schmidt, and B.~Haasdonk}, {\em Model order reduction
  approaches for infinite horizon optimal control problems via the hjb
  equation}, in Model Reduction of Parametrized Systems, P.~Benner,
  M.~Ohlberger, A.~Patera, G.~Rozza, and K.~Urban, eds., Springer International
  Publishing, Cham, 2017, pp.~333--347.

\bibitem{A05}
{\sc A.~Antoulas}, {\em Approximation of Large--Scale Dynamical Systems},
  {SIAM} Publications, Philadelphia, PA, 2005.

\bibitem{AK01}
{\sc J.~Atwell and B.~King}, {\em Proper orthogonal decomposition for reduced
  basis feedback controllers for parabolic equations}, Mathematical and
  Computer Modelling, 33 (2001), pp.~1--19.

\bibitem{Barkouki2015}
{\sc H.~Barkouki, A.~H. Bentbib, and K.~Jbilou}, {\em An adaptive rational
  method{L}anczos-type algorithm for model reduction of large scale dynamical
  systems}, Journal of Scientific Computing, 67 (2015), p.~221–236.

\bibitem{GB17}
{\sc C.~Beattie and S.~Gugercin}, {\em Model reduction by rational
  interpolation}, in Model Reduction and Approximation: Theory and Algorithms,
  SIAM, Philadelphia, 2017.

\bibitem{Benner.Bujanovic.16}
{\sc P.~Benner and Z.~Bujanovi{\'c}}, {\em On the solution of large-scale
  algebraic {Riccati} equations by using low-dimensional invariant subspaces},
  Linear Algebra and Its Applications, 488 (2016), pp.~430--459.

\bibitem{BCOW.17}
{\sc P.~Benner, A.~Cohen, M.~Ohlberger, and K.~Willcox}, eds., {\em Model
  reduction and approximation: theory and Algorithms}, Computational Science \&
  Engineering, SIAM, PA, 2017.

\bibitem{Benner.Li.Penzl.08}
{\sc P.~Benner, J.-R. Li, and T.~Penzl}, {\em Numerical solution of large-scale
  {L}yapunov equations, {R}iccati equations, and linear-quadratic optimal
  control problems}, Num. Lin. Alg. with Appl., 15 (2008), pp.~1--23.

\bibitem{BS10}
{\sc P.~Benner and J.~Saak}, {\em A {Galerkin-Newton-ADI} method for solving
  large-scale algebraic {R}iccati equations.}, Tech. Rep. 1253, DF Priority
  program, 2010.

\bibitem{Benner.Saak.survey13}
{\sc P.~Benner and J.~Saak}, {\em Numerical solution of large and sparse
  continuous time algebraic matrix {R}iccati and {L}yapunov equations: A state
  of the art survey}, GAMM-Mitt.,  (2013), pp.~32--52.

\bibitem{bG15}
{\sc J.~Borggaard and S.~Gugercin}, {\em Model reduction for {DAE}s with an
  application to flow control}, Active Flow and Combustion Control 2014, R.
  King editors, Springer-Verlag, Notes on Numerical Fluid Mechanics and
  Multidisciplinary Design,, 127 (2015), pp.~381--396.

\bibitem{Druskin.Simoncini.11}
{\sc V.~Druskin and V.~Simoncini}, {\em Adaptive rational {Krylov} subspaces
  for large-scale dynamical systems}, Systems and Control Letters, 60 (2011),
  pp.~546--560.

\bibitem{FF14}
{\sc M.~Falcone and R.~Ferretti}, {\em Semi-{L}agrangian Approximation Schemes
  for Linear and {H}amilton-{J}acobi equations}, SIAM, 2014.

\bibitem{GP11}
{\sc L.~Gr{\"u}ne and J.~Pannek}, {\em Nonlinear model predictive control.
  Theory and algorithms}, Springer Verlag, 2011.

\bibitem{Gugercin2008a}
{\sc S.~Gugercin, A.~C. Antoulas, and C.~Beattie}, {\em {${\cal H}_2$ model
  reduction for large-scale linear dynamical systems}}, SIAM J. Matrix Anal.
  Appl., 30 (2008), pp.~609--638.

\bibitem{GutknechtApr.1992}
{\sc M.~H. Gutknecht}, {\em {A} completed theory of the unsymmetric {L}anczos
  process and related algorithms. {I}}, SIAM J. Matrix Anal. Appl., 13 (Apr.
  1992), pp.~594--639.

\bibitem{Heyouni.Jbilou.09}
{\sc M.~Heyouni and K.~Jbilou}, {\em An extended {Block Krylov method} for
  large-scale continuous-time {algebraic Riccati equations}}, ETNA, 33
  (2008-2009), pp.~53--62.

\bibitem{HPUU09}
{\sc M.~Hinze, R.~Pinnau, M.~Ulbrich, and S.~Ulbrich}, {\em Optimization with
  {PDE} Constraints. Mathematical Modelling: Theory and Applications}, Springer
  Verlag, New York, 2009.

\bibitem{Jaimoukha.Kasenally.95}
{\sc I.~M. Jaimoukha and E.~M. Kasenally}, {\em Oblique projection methods for
  large scale model reduction}, SIAM J. Matrix Anal. Appl., 16 (1995),
  pp.~602--627.

\bibitem{JaimoukhaFeb.1994}
{\sc I.~M. Jaimoukha and E.~M. Kasenally}, {\em {Krylov subspace methods for
  solving large {Lyapunov} equations}}, SIAM J. Numer. Anal., 31 (Feb. 1994),
  pp.~227--251.

\bibitem{Jbilou_03}
{\sc K.~Jbilou}, {\em {Block Krylov subspace methods for large algebraic
  Riccati equations}}, {Numerical Algorithms}, 34 (2003), pp.~339--353.

\bibitem{K68}
{\sc D.~Kleinman}, {\em On an iterative technique for riccati equation
  computations}, IEEE Trans. Automat. Control, 13 (1968), pp.~114--115.

\bibitem{KS16}
{\sc B.~Kramer and J.~Singler}, {\em A pod projection method for large-scale
  algebraic riccati equations.}, Numer. Algebra Control Optim., 4 (2016),
  pp.~413--435.

\bibitem{KV08}
{\sc K.~Kunisch and S.~Volkwein}, {\em Proper orthogonal decomposition for
  optimality systems}, ESAIM: Mathemematical Modelling and Numerical Analysis,
  42 (2008), pp.~1--23.

\bibitem{Lancaster.Rodman.95}
{\sc P.~Lancaster and L.~Rodman}, {\em Algebraic {R}iccati equations}, Oxford
  Univ. Press, 1995.

\bibitem{Lin.Simoncini.12}
{\sc Y.~Lin and V.~Simoncini}, {\em {Minimal residual methods for large scale
  Lyapunov equations}}, tech. rep., Dipartimento di Matematica, Universit\`a di
  Bologna, 2012.

\bibitem{Lin.Simoncini.15}
{\sc Y.~Lin and V.~Simoncini}, {\em A new subspace iteration method for the
  algebraic {R}iccati equation}, Numerical Linear Algebra w/Appl., 22 (2015),
  pp.~26--47.

\bibitem{Q12}
{\sc A.~Quarteroni}, {\em Modellistica numerica per problemi differenziali},
  Springer, 2012.

\bibitem{SH17}
{\sc A.~Schmidt and B.~Haasdonk}, {\em Reduced basis approximation of large
  scale parametric algebraic {Riccati} equations}, ESAIM: Control, Optimisation
  and Calculus of Variations,  (2017).

\bibitem{Simoncini2016}
{\sc V.~Simoncini}, {\em Analysis of the rational {K}rylov subspace projection
  method for large-scale algebraic {R}iccati equations}, SIAM J. Matrix Anal.
  Appl., 37 (2016), pp.~1655--1674.

\bibitem{Simoncini.16}
\leavevmode\vrule height 2pt depth -1.6pt width 23pt, {\em On the extended
  {K}rylov subspace method for the algebraic {R}iccati equation}, January 2016.
\newblock In preparation.

\bibitem{SSM14}
{\sc V.~Simoncini, D.~B. Szyld, and M.~Monslave}, {\em On two numerical methods
  for the solution of large-scale algebraic riccati equations}, IMA, J. Numer.
  Anal., 34 (2014), pp.~904--920.

\bibitem{Sir87}
{\sc L.~Sirovich}, {\em Turbulence and the dynamics of coherent structures.
  parts i-ii}, Quarterly of Applied Mathematics,, XVL (1987), pp.~561--590.

\bibitem{Vol11}
{\sc S.~Volkwein}, {\em Model Reduction using Proper Orthogonal Decomposition},
  Lecture notes, University of Konstanz, 2011.

\end{thebibliography}

\end{document}